\documentclass[12pt]{article}

         \input xy
         \xyoption{all}

         \usepackage{amsfonts}
         \usepackage{latexsym}

         \newtheorem{thm}{Theorem}[section]

         \newtheorem{lem}[thm]   {Lemma}
         \newtheorem{cor}[thm]   {Corollary}
         \newtheorem{rem}[thm]   {Remark}
         \newtheorem{defn}[thm]  {Definition}
         \newtheorem{ex}[thm]    {Example}
         \newtheorem{prop}[thm]  {Proposition}

         \newcommand{\term}[1]   {{\bf #1}}

         \newenvironment{proof}  {\par\noindent{\bf Proof}\ }
                                  {\hfill$\Box$\par\medskip}

         \newcommand{\maprt}[1]
{\smash{\mathop{\longrightarrow}\limits^{#1}}}
         
         \newcommand{\inclds}     {\hookrightarrow}

         \newcommand{\A}        {\mathcal{A}}

         \newcommand{\F}        {\mathcal{F}}

         \newcommand{\E}        {\mathcal{E}}
\newcommand{\SSS}        {\mathcal{S}}

         \newcommand{\CC}    {\Bbb{C}}

         \newcommand{\id}    {\mathrm{id}}

         \newcommand{\twdl}  {\widetilde}

         \newcommand{\sseq}  {\subseteq}

         \newcommand{\cp}    {\CC\mathrm{P}}

         \newcommand{\cat}   {\mathrm{cat}}

\newcommand{\sct}{\mathrm{secat}}

\newcommand{\dash}{\vrule height3pt depth-2.3pt width.1cm\hskip.05cm}

\def\shortlabeldecomp#1#2#3#4#5#6#7#8{
\xymatrix{
{#1}_0\ar[r]^{{#7}_0} \ar@{=}[d] & {#1}_1\ar[r]^{{#7}_1}&
\cdots \ar[r]^(.42){{#7}_{{#4}-2}}
  & {#1}_{#4 -1 }\ar[r]^(.52){{#7}_{{#4}-1}} & {#1}_{#4}
\ar@<-.4ex>[d]_{{#2}_{#4}}\\
#1\ar[rrrr]^{#2} &&&& #3 \ar@<-.4ex>[u]_{{#8}} \\ }
}

         \title{The Sectional Category of a Map}
         \author{M. Arkowitz and J. Strom}

         \date{}   

\begin{document}

         \maketitle

\abstract{\noindent
We study a generalization of the Svarc genus of a fiber map.  For an
arbitrary collection $\E$ of spaces and a map $f:X\to Y$, we define a
numerical invariant, the $\E$-sectional category of $f$,
in terms of open
covers of $Y$.  We obtain several basic properties of
$\E$-sectional category, including
those dealing with homotopy domination and homotopy pushouts.  We then 
give
three simple axioms which characterize the
$\E$-sectional category.  In the final section we
obtain inequalities for the
$\E$-sectional category of a composition and inequalities
relating the
$\E$-sectional category to the Fadell-Hussein category of a map and the
Clapp-Puppe category of a map.
}

\bigskip


\smallskip



\bigskip


\section{Introduction}

The sectional category of  a fiber map $f:X\to Y$, denoted $\sct (f)$,
is one less than the number of sets in the smallest
open cover of $Y$ such that $f$ admits a cross-section over each member
of
the cover.
This simple and natural numerical invariant of fiber maps was
developed and studied extensively by Svarc \cite{Sv} who
called it the genus of $f$.
Subsequently
Fet \cite{Fe} and
Berstein-Ganea \cite{BG} extended the definition to arbitrary maps and
related it to
the Lusternik-Schnirelmann category.  There have been many applications
of
sectional category to questions of classification of bundles, embeddings
of
spaces and existence of regular maps \cite{Sv} as well as
applications
outside of
algebraic and differential topology \cite{Da,FW}.  However, since
Svarc's papers, the
actual study of sectional category has been sporadic and has often
appeared as
subsidiary results within a larger work
\cite[\S 2]{BG},
  \cite[\S 8]{Ja},
  \cite[\S 4]{CP},
\cite[9.3]{CLOT}.
(An exception to this is the paper of Stanley \cite{St} which deals
with
sectional
category in the context of rational homotopy theory.)

Recently there has been considerable interest in the study of numerical
homotopy invariants of spaces and of maps.
Classically the Lusternik-Schnirelmann category and cone length of a 
space
have been studied in  \cite{LS,Fo,Ga,Co1} and, more
recently, the cone length \cite{Mar}, Clapp-Puppe category
\cite{CP} and
Fadell-Husseini category of a map \cite{FH,Co2} have been
investigated.

Together with D. Stanley, we gave a unified axiomatic development of
many of these invariants in the paper \cite{ASS}.
For the category of spaces and
maps and a fixed collection $\E$ of objects, these axioms were of two 
basic
types, namely, those dealing with numerical functions relative to
$\E$ defined on the objects of the category and those dealing with
numerical functions relative to $\E$ defined on the morphisms of the
category.  By specializing the collection $\E$, we obtained some of the
previously  studied invariants as well as several new invariants.  For
instance, by setting
$\E= \{\mathrm{all\ spaces}\}$,
our axioms for  morphisms yield  the Fadell-Husseini
category of a map. One invariant that was not dealt with in \cite{ASS}
was the sectional category.  Whereas the various versions of category
and cone length have been defined in numerous homotopy-invariant ways,
this is not the case for sectional category.  The fact that sectional
category does
not fit nicely into the general framework that so neatly
encapsulates the category and cone length of spaces and maps may
account for its marginal status in the decades since its introduction.

We begin in Section 3 by defining a generalization of the sectional
category
of a map
with respect to a collection $\E$ of spaces using open covers of
the target of the map.  This is a straightforward extension of the
classical
definition.  We derive several simple basic properties of this
invariant. In
particular, when
$\E =\{\mathrm{all\ spaces}\}$, we see that
$\E\dash\sct(f) =  \sct(f)$.
This treatment of sectional category leads to new invariants obtained by
varying the  collection $\E$ of spaces.

In Section 4 we bring sectional category of maps in line with the
other invariants studied in \cite{ASS}.  This is done by considering
maps as {\it objects} in the category
whose objects are maps of spaces and whose morphisms are given by
commutative
squares.
In this category, we apply the  axiomatic approach for invariants
of objects,
   relative to  the collection   of all maps with sections which
factor through a space in $\E$.
The unique numerical invariant obtained from
the axioms is then proved to be $\E\dash\sct$.  This is the content of
Theorem
\ref{thm:main}.

Another case of this axiomatic approach is obtained
by working in  the category  of maps, relative to the
the collection  of all maps which factor through
some space in $\E$.
The unique invariant obtained from the axioms
can then be shown to be
the $\E\dash$Clapp-Puppe category of a map (see Remark 
\ref{rem:duality}).  We
hope to
return to this in the future.

All invariants that fall into our axiomatic scheme will of course
share certain basic properties that follow formally from the axioms.
A major interest, however, is in the new questions
which arise
regarding $\E\dash$sectional category.
Many of these are considered in Section 5 and take the form of
inequalities.
There we concentrate on the following:
(1)  How is $\E\dash\sct(f)$ related to the domain and the target of $f$?
(2)  How does $\E\dash$sectional category behave  with respect to
       composition of maps?
(3) What is the relation between the $\E\dash$sectional category,
       the Clapp-Puppe category and the Fadell-Husseini category?

\section{Preliminaries}\label{section:preliminaries}

In this section we establish our notation and recall some
basic definitions and results which we shall use.  All spaces are to have
the homotopy type of connected CW-complexes.  We do not assume that
spaces have base points and hence maps are not base point preserving.
For spaces $X$ and $Y$, $X\equiv Y$ denotes that $X$ and $Y$ have the
same homotopy type.  We let $*$ denote a space with one point and we
also write $*:X\to Y$ for any constant map from $X$ to $Y$.  We denote
by $\id_X$ or id the identity map of $X$.  If $f,g:X\to Y$ are two maps,
then $f\simeq g$ signifies that $f$ and $g$ are homotopic.  Given maps
$f:X\to Y$ and $g:Y\to X$, if $fg \simeq \id _Y$, we say that $g$ is a
{\bf
section} of $f$ or that $Y$ is a {\bf retract} of $X$.  For maps $f:X\to 
Y$
and $f':X'\to Y'$,
if there is a homotopy
commutative diagram
$$
\xymatrix{
X\ar[rr]^i\ar[d]^f && X'\ar[d]^{f'} \ar[rr]^r && X\ar[d]^f\\
Y\ar[rr]^j  && Y'  \ar[rr]^s && Y \\ }
$$
such that $ri \simeq \id$ and $sj \simeq \id$,
then we say that $f'$ {\bf dominates} $f$ or
that $f$ is a {\bf retract} of $f'$.
If, in addition,
$ir\simeq\id$ and $js \simeq \id$, then $f$ and $f'$ are called
{\bf equivalent}, and we write $f \equiv f'$.

We will also use certain basic constructions in homotopy
theory. The \term{pushout} of a diagram of the form
$$
\xymatrix@1{
C && A\ar[ll]_g \ar[rr]^f && B }
$$
is the quotient space of $B\cup C$ by  the
equivalence relation which sets $f(a)$
equivalent to $g(a)$, for every $a\in A$.
The \term{homotopy pushout} $H$ of the diagram
is the
quotient space of $B\cup (A\times [0,1])\cup C$ by  the
equivalence relation which sets $(a,0)$
equivalent to $f(a)$ and $(a,1)$
equivalent to $g(a)$, for every $a\in A$.
The pushout and homotopy pushout
constuctions
are clearly functors from the category of given diagrams to
the category of spaces.

There are two homotopy pushouts
of special interest.  The first is the homotopy pushout of
$$
\xymatrix@1{{* } && A\ar[ll]\ar[rr]^f && B }
$$
which is the \term{mapping cone} of $f$ and is
denoted $B\cup CA$.  The second is the homotopy
pushout of
$$
\xymatrix@1{
A && A\ar[ll]_{\id _A} \ar[rr]^f && B }
$$
which is the \term{mapping cylinder} of $f$ and is
denoted $M_f$.  We note for later use that the homotopy
pushout of
$$
\xymatrix@1{
C && A\ar[ll]_g \ar[rr]^f && B }
$$
is homeomorphic to the pushout of the associated mapping cylinder diagram
$$
\xymatrix@1{
M_g && A\ar[ll]_{} \ar[rr]^{} &&M_f }.
$$

We next consider two versions of
the category of a map.  Given $f:X\to Y$, we first
define the (reduced) \term{Fadell-Husseini category}
of $f$, denoted $\cat_{\mathrm{FH}}(f)$ \cite{FH}.  We note
that $f$ is equivalent to the inclusion of $X$ into $M_f$,
and so we regard $f$ as an inclusion.  Then $\cat_{\mathrm{FH}}(f)$
is the smallest $n$ such that there is an
open cover  $\{U_0,U_1,\ldots ,U_n\}$ of $Y$ with the following
properties:
(1) if $j_i:U_i\to Y$ is the inclusion, then $j_i \simeq *$
for $i= 1,..., n$,  (2)
$X\subseteq U_0$
and  (3) there is a map
$r:U_0\to X$ and a homotopy of pairs $j_0\simeq r:(U_0, X)\to
(Y,X)$.
It follows that $\cat_{\mathrm{FH}}(*\to Y)$ is just the (reduced)
Lusternik-Schnirelmann category $\cat (Y)$
of $Y$.

By a {\bf collection} $\E$ we mean any collection of
spaces which contains the one point space $*$.
The second notion of category of a map $f:X\to Y$ that we
consider is the (reduced) $\E\dash$\term{Clapp-Puppe category} of
$f$, denoted
$\E\dash\cat_{\mathrm{CP}}(f)$ \cite{CP}.  This is the smallest 
non-negative
integer $n$ with the
following properties: (1) there exists an open cover
$\{U_0,U_1,\ldots ,U_n\}$ of
$X$,  (2) there exists spaces $E_i \in \E$ and maps $u_i: U_i\to E_i$
and $v_i:E_i\to Y$, for $i=0,1,\ldots ,n$  and (3) $f|U_i \simeq
v_iu_i$  for each $i$.  For a space $X$, the
$\E\dash$Clapp-Puppe category
of
$X$ is defined by $\E\dash\cat_{\mathrm{CP}}(X) =
\E\dash\cat_{\mathrm{CP}}({\mathrm{id}}_X)$.
If $\E = \{*\}$, the collection consisting of a one point
space, then $\E\dash\cat_{\mathrm{CP}}(f)$ is the category of the
map $f$, as discussed by Berstein-Ganea \cite{BG} and others, which
we will denote by $\cat_{\mathrm{BG}}(f)$.  We note that
$\cat_{\mathrm{BG}}({\mathrm {id}}_X) = \cat (X)$.

\section{Definition and Basic Properties}

In this section we define the sectional category
of a map relative to a collection $\E$ of spaces.  We then establish
basic properties of this invariant.  Of particular importance are the
Domination Proposition (Prop.\,\ref{prop:domination}) and the Homotopy
Pushout Theorem (Thm.\,\ref{thm:po}).

\begin{defn}\label{defn:defs}
\begin{enumerate}
{\em
\item
If $\E $ is a collection, $f:X \to Y$ a map and
$i:U \inclds Y$ the inclusion map, then $U$ is
\term{$\E$-section-categorical}
(with respect to $f$) if there is space
$E \in \E$ and maps $u:U\to E$ and $v:E \to X$
such that the following diagram is homotopy-commutative:
$$
\xymatrix{
E \ar[rr]^{v}               &&   X \ar[d]^{f}  \\
U \ar[rr]^{i} \ar[u]^{u}    &&   Y.            \\  }
$$
The map $vu:U \to X$ is called an {\bf $\E\dash$section of $f$ over $U$}.
\item
For a map $f:X\to Y$, the {\bf (reduced)
$\E$-sectional category } of $f$,
written $\E\dash\sct (f)$, is   the
smallest   integer $n$ such that there exists an
open cover of $Y$ by $n+1$ subsets, each of which is
$\E$-section-categorical.
If no such integer exists, then $\E\dash\sct (f) = \infty $.  }
\end{enumerate}
\end{defn}
If $\E$ is the collection of  all spaces, then we write
$\sct(f)$ for $\E\dash\sct (f)$ and note that $\sct(f)$
is just the sectional category of $f$ as defined in \cite[Def.\thinspace
2.1]{BG}.

\medskip

The following result lists a number of useful
properties of $\E\dash\sct$.  The proofs are all straighforward,
and we omit them.

\begin{prop}\label{prop:simple}
Let $\E$ be a collection.
\begin{enumerate}
\item
If $f\simeq f'$, then
$\E\dash\sct (f) =\E\dash\sct (f'$).
\item
If $\E = \{ * \} $, then $\E\dash\sct (f)  = \cat(Y)$,
the category of $Y$.
\item
If $\E$ and $\F$ are two collections and $\F \subseteq \E$, then
$\E\dash\sct (f) \le \F\dash\sct (f)$.
\item
If $*:X\to Y$ is a constant map, then
$\E\dash\sct (*) = \cat (Y)$.
\end{enumerate}
\end{prop}

\noindent
Assertions (2) and (4) show that the Lusternik-Schnirelmann category of
$Y$ can be obtained from the $\E$-sectional
category of $f$ in two ways: either by making $\E$ trivial or by
making $f$ trivial.
Since $\{ *\}\subseteq \E \subseteq  \{ \mathrm{all\ spaces}\}$,
we obtain from (3) the following.

\begin{cor}\label{cor:id}
For any map $f:X\to Y$ and collection $\E $,
$$
\sct(f) \leq \E\dash\sct(f) \leq \cat(Y).
$$
\end{cor}

Next we consider some basic properties of the identity map.

\begin{prop}\label{prop:id}
Let $\E$ be a collection and $X$ be a space.
\begin{enumerate}
\item
$\E\dash\sct(\id_X) = \E\dash\cat_{\mathrm{CP}}(X)$.
\item
$X$ is a retract of a space in $\E$ if and only if $\E\dash\sct
(\id_X) = 0$.  In particular, $\sct(\id_X) = 0$ for every space $X$.
\end{enumerate}
\end{prop}

The following examples show that the inequalities in Proposition
\ref{prop:simple} and Corollary \ref{cor:id} can be strict.

\begin{ex}
{\em
\begin{enumerate}
\item
For the inequality in \ref{prop:simple}, let $\E$ be the collection of
all
spaces and $\F$ a collection which does not contain a space
having $X$ as a retract.
Then
$$
\E\dash\sct(\id_X) = 0< \F\dash\sct(\id_X).
$$
For a specific example, let $\F = \{ *\}$ and take any $X\not\equiv *$.
\item
For the inequality in \ref{cor:id}, consider
$f:X\to Y$ and  the inclusion $j: Y \to Y\cup CX$.  Then $\sct(j) \leq 
1$ by Lemma
\ref{lem:mc} (below). If $\cat(Y\cup CX)> 1$, then
$$
\sct(j) \leq  1  < \cat(Y\cup CX).
$$
This is the case, for example, when
$f: S^{2n+1} \to \cp^n$ is the Hopf map.
\end{enumerate}
}
\end{ex}

We next establish a basic property of $\E$-sectional
category.  In the next section this will be shown to be
one of  three properties which characterize $\E\dash\sct$.

\begin{prop}\label{prop:domination} {\em {(Domination)}}
If $f:X \to Y$ is dominated by $f':X'\to Y'$, then
$$
\E\dash\sct (f) \le \E\dash\sct (f').
$$
\end{prop}

\begin{proof}
We are given $j:Y\to Y'$ and $r:X'\to X$ as in the definition
of domination in \S 2.
Let $\E\dash\sct (f') =n$ with
$\E$-section-categorical  cover  $\{ U'_0,U'_1,\ldots ,U'_n\} $
of $Y'$ and maps
$$
\xymatrix@1{
U'_i\ar[rr]^{u_i} && E_i\ar[rr]^{v_i} && X'
}
$$
such that $f'v_iu_i \simeq j_i:U'_i\to Y'$, where $E_i \in \E$
and $j_i$ is the inclusion.
If $U_i = j^{-1}(U'_i)$, then $\{ U_0,U_1,\ldots ,U_n\} $
is an open cover of $Y$
and
$$
\xymatrix@1{
U_i     \ar[r]^{j|U_i} &
U'_i    \ar[r]^{u_i}   &
E_i     \ar[r]^{v_i}   &
X'      \ar[r]^{r}     &    X
}
$$
is the desired $\E\dash$section of $f$ over $U_i$.
\end{proof}

An immediate corollary of Proposition \ref{prop:domination}
is that $\E\dash\sct$ is an invariant of homotopy equivalence of
maps.

\begin{cor}  \label{cor:equiv}
If $f \equiv f'$, then $\E\dash\sct (f) =\E\dash\sct (f')$.
\end{cor}

It is well-known that every map is homotopy equivalent
to a fiber map \cite[p.\,48]{May}.
Corollary \ref{cor:equiv} implies that the
$\E$-sectional category of an arbitrary map is equal to
the $\E$-sectional category
of the equivalent fiber map.  Therefore Svarc's
definition of sectional category \cite{Sv}, which
applies only to fiber maps, is equivalent to
Definition \ref{defn:defs} in the special case
$\E= \{ \mathrm{all\ spaces} \}$.

\medskip

We next prove a result about the $\E\dash$sectional category of the
maps of one homotopy pushout into another.
To establish notation, let
$$
\xymatrix{
C\ar[d]^{c} && A \ar[ll]_{g} \ar[rr]^{f} \ar[d]^{a} && B\ar[d]^{b}
\\
C' && A' \ar[ll]_{g'}\ar[rr]^{f'} && B'
\\  }
$$
be a commutative diagram, let $D$ and $D'$ be the homotopy
pushouts of the top and bottom rows, respectively, and let $d:D\to D'$
be the induced map.  We begin with a lemma.

\begin{lem}\label{lem:cover}
If $\E\dash\sct (b) = n$,
then there exists $\E$-section-categorical open
sets $N_0,N_1,\ldots N_n$
of $D'$ with respect to $d$
such that $M_{f'}\subseteq \bigcup N_i $.
\end{lem}

\begin{proof}
Let $\twdl d: M_f \to M_{f'}$ be the map $d$ with restricted domain and
target.  Then $\twdl d \equiv b$, so $\E\dash\sct (\twdl d) = \E\dash\sct
(b) = n$
by Corollary \ref{cor:equiv}.
To complete the proof,   observe that if $U\sseq  M_{f'}$ is
$\E$-section-categorical
with respect to $\twdl d$, then it is also $\E$-section-categorical
with respect to
$d$.
\end{proof}

Now we can prove the Homotopy Pushout Theorem, which is another
basic property which we use in \S\ref{section:axioms} to characterize
$\E\dash $sectional category.

\begin{thm}\label{thm:po} {\em (Homotopy Pushout)}
With the notation above,
$$
\E\dash\sct (d) \le \E\dash\sct (b) + \E\dash\sct (c) +1 .
$$
\end{thm}

\begin{proof}
Let $\E\dash\sct (b) = n $ and $\E\dash\sct (c) = m $.
Let $\{ U_0, U_1,\ldots , U_n\} $ be a
minimal open $\E$-section-categorical cover  of $B'$ with respect to $b$.
Then there are $\E$-section-categorical
open sets
$N_0,N_1, \ldots ,N_n$ of $D'$ with respect to $d$ which
cover $M_{f'}$ by Lemma \ref{lem:cover}.  Similarly,
if $\{ V_0, V_1,\ldots , V_m\} $ is a
minimal open $\E$-section-categorical cover  of $C'$ with respect to $c$,
then there are $\E$-section-categorical open sets
$M_0,M_1, \ldots ,M_m$ of $D'$
with respect to $d$ which cover
$M_{g'}$.  It follows that $\{ N_0, \ldots ,N_n, M_0, \ldots ,M_m\}$ is
an
$\E$-section-categorical cover of $D'$
with respect to $d$.  Thus $\E\dash\sct (d) \le \E\dash\sct (b) +
\E\dash\sct (c) +1 $.
\end{proof}

We conclude this section with a simple application of Theorem
\ref{thm:po} to the maps in a homotopy pushout square.

\begin{cor}
If
$$
\xymatrix{
A\ar[rr]^{f}   \ar[d]_{g}     &&   B \ar[d]^{r}  \\
C \ar[rr]^{s}     &&   D   \\  }
$$
is a homotopy pushout square, then
$$
\E\dash\sct (r) \le \E\dash\sct (g)
+\E\dash\cat_{\mathrm{CP}}(B)
+1.
$$
\end{cor}

\begin{proof}
The map of homotopy pushouts obtained from
the commutative diagram
$$
\xymatrix{
A\ar[d]^{g}   && A  \ar@{=}[ll]  \ar[rr]^{f} \ar@{=}[d]  && B\ar@{=}[d]
\\
C && A \ar[ll]_{g}\ar[rr]^{f} && B
\\  }
$$
is equivalent to $r:B \to D$.  Now apply Theorem \ref{thm:po},
using Proposition \ref{prop:id}(1).
\end{proof}

\section{Axioms}\label{section:axioms}

In this section we characterize the $\E\dash$sectional
category by simple axioms.  For reasons given in Remark
\ref{rem:duality} we state our axioms in
greater generality than is needed for sectional category.

\begin{defn}\label{defn:axi}
{\em We denote by $\SSS$ a non-empty collection of maps.
An \term{$\SSS$-length function} is a function
$\gamma = \gamma_{\SSS}$ which assigns to every map $f$ an integer
$0\le \gamma (f) \le \infty $ such that
\begin{enumerate}
\item
If $f\in \SSS$, then $\gamma (f) = 0$.
\item
Let
$$
\xymatrix{
C\ar[d]^{c} && A \ar[ll]_{g} \ar[rr]^{f} \ar[d]^{a} && B\ar[d]^{b}
\\
C' && A' \ar[ll]_{g'}\ar[rr]^{f'} && B'
\\  }
$$
be a commutative diagram with induced map $d:D\to D'$
of the homotopy pushout of the first row into the homotopy
pushout of the second row.
If $c\in \SSS $, then $\gamma (d) \le \gamma (b) +1$.
\item If $f$ is dominated by $f'$, then
$\gamma (f) \le \gamma (f')$.
\end{enumerate}
We call (1)--(3) the {\bf axioms for $\SSS$-length functions}.  }
\end{defn}

It is an immediate consequence of Axiom (3) that if $f\equiv f'$, then 
$\gamma (f) = \gamma (f')$.
Then, since homotopic maps are equivalent (\S 2), it follows that
$f\simeq f'$ implies that $\gamma (f) = \gamma (f')$.

\begin{rem}\label{rem:categorical}
{\em
These axioms are analogous to the
axioms for the $\A$-category of a {\it space}
given in \cite[Prop.\thinspace 5.6(2)]{ASS}.
The following comments are made to elucidate the analogy.
In Definition \ref{defn:axi}
we define the $\SSS$-length function on the objects in the category
of maps of spaces.
In \cite[p.\,26]{ASS} an $\SSS$-length function was defined on the 
objects
in the category of spaces in the case where $\SSS $ is the collection
of contractible spaces, though we could have taken an arbitrary
collection $\SSS$ in the definition.
This is why the analog of Axiom (3) in \cite{ASS} deals with mapping
cone sequences
instead of homotopy pushouts.  In addition, the axioms in 
\cite[p.\,26]{ASS}
were made with respect to
an arbitrary collection $\A$ of spaces.  It is possible to incorporate an
arbitrary collection $\A$ of maps into
Definition \ref{defn:axi}, but we have not given the definition in this 
generality.}
\end{rem}

\noindent  We are led to the following definition.

\begin{defn} \label{defn:maxi}
{\em
The \term{maximal $\SSS$-length function} $\Gamma_{\SSS}$  is defined by
$$
\Gamma_\SSS (f) = \max
\{
\gamma_\SSS (f)\: | \:\mathrm{for\ every}
\: \SSS\dash\mathrm{length\ function}
\: \gamma_\SSS
\}   .
$$}
\end{defn}

\noindent

Note that $\Gamma_\SSS$ satisfies the axioms  for $\SSS$-length
functions.

\medskip

We wish to show that $\E\dash\sct$ equals
$\Gamma_\SSS$  for some $\SSS$.  For this we need to   choose
a collection $\SSS$ of maps that is closely related to  $\E$.

\begin{defn}
{\em  For a collection of spaces $\E$, let
$$
\SSS(\E)
=
\{ f:X\to Y\, |\,
\mathrm{there\ is\ an}\ \E\dash\mathrm{section\ of}\  f\
\mathrm{over}\  Y\}.
$$
}
\end{defn}

\noindent
The main result of this section is that $\Gamma_{\SSS(\E)} =
\E\dash\sct$.
In order to prove this, we need to show that we can replace open covers
with covers
by simplicial subcomplexes in the definition of $\E\dash$sectional
category.
This is crucial in dealing with maps that are defined by homotopy
pushouts.
We require  the following lemma.

\begin{lem}\label{lem:sub}
Let $f:X\to Y$ be a map and $B\subseteq Y$ an
$\E$-section-categorical subset.   If
$A\subseteq B$, then $A$ is $\E$-section-categorical.
If $B\subseteq C$ is a deformation
retract, then $C$ is $\E$-section-categorical.
\end{lem}

The proof is elementary, and we omit it.

        \begin{prop}\label{prop:scx}
        Let $Y$ be a  simplicial complex and $f: X\to Y$ a map.
        Then $\E\dash\sct (f) \le n \iff $ there is a simplicial structure
        on
        $Y$ such that $Y$ can be covered by $n+1$  
$\E$-section-categorical
        subcomplexes $L_0,L_1,\ldots ,L_n$.
        \end{prop}

\begin{proof}
$(\Leftarrow)$\quad If $N_i$ is the second regular
neighborhood of $L_i$ in $Y$, then $L_i$
is a deformation retract of the closure $\overline N_i$ 
\cite[p.\,72]{ES}.
Thus $\overline N_i$ is $\E$-section-categorical,
and hence so is the open set $N_i$ by Lemma \ref{lem:sub}.

$(\Rightarrow)$\quad   Consider an $\E$-section-categorical
cover $\{ U_0,U_1, \ldots ,U_n \}$ of $Y$.  It is known that
there is a simplicial structure on $Y$ with subcomplexes
$L_0,\ldots , L_n$
which cover $Y$ such that $L_i\sseq U_i$ for each $i$
\cite[Lem.\,2.3]{GH}.
The result follows from Lemma \ref{lem:sub}.
\end{proof}

\begin{thm}\label{thm:main}
If $\E$ is a collection of spaces and
$f:X\to Y$ is a map, then
$$
\Gamma _{\SSS (\E )}(f) = \E\dash\sct (f).
$$
\end{thm}

\begin{proof}
In the proof we write $\Gamma$ for $\Gamma_{\SSS(\E )}$.
To show $\E\dash\sct (f) \le \Gamma  (f)$,
it suffices to show that $\E\dash\sct $ satisfies the axioms of
Definition
\ref{defn:axi} for an $\SSS(\E) $--length function.  But this follows 
from
\S 3.

Next we show $\Gamma (f) \le \E\dash\sct (f)$.
We suppose that $\E\dash\sct (f) < \infty $ and
prove the result by induction on $\E\dash\sct (f)$.
If $\E\dash\sct (f) = 0$, then $f \in \SSS(\E)$,  and so $\Gamma (f) =0$.
Suppose next that the inequality holds for all maps $g$ with
$\E\dash\sct (g) < n $
and let $f:X\to Y$  with $\E\dash\sct (f)= n$.
Now $X$ and $Y$ have the homotopy type of simplicial complexes
\cite[Thm.\,2]{Mi},
so we may assume that $X$ and $Y$ are simplicial complexes.
By the Simplicial Approximation Theorem,
we can take  $f$  to be  a simplicial map.
By Proposition \ref{prop:scx}, there is an $\E$-section-categorical
cover $\{L_0,L_1,\ldots ,L_n\}$ of $Y$ by subcomplexes with respect to
$f$. We set $K_i = f^{-1}(L_i)$ so that
$\{K_0,K_1,\ldots ,K_n\}$ is a cover of $X$ by subcomplexes.

Thus the following diagram is commutative
$$
\xymatrix{
K_n\ar[d]^{f_1} && K_n\cap (K_0\cup \dots \cup K_{n-1}) \ar[ll]_{}
\ar[rr]^{} \ar[d]^{f_2} &&
K_0\cup \dots \cup K_{n-1}\ar[d]^{f_3}
\\
L_n && L_n\cap (L_0\cup \dots \cup L_{n-1}) \ar[ll]_{}\ar[rr]^{ } &&
L_0\cup \dots \cup L_{n-1},
\\  }
$$
where each $f_i$ is induced by $f$.  Clearly $f_1 \in \SSS(\E )$.  
Furthermore,
$\E\dash\sct (f_3) \le n-1$, and hence $\Gamma (f_3)
\le n-1$ by our inductive hypothesis.  But the pushout of the rows of the
diagram
yields the map $f:X\to Y$.  However
this is homotopy equivalent to the homotopy
pushout since all horizontal maps are
simplicial inclusions and hence cofibrations \cite[p.\,78]{May}.
By Definition \ref{defn:axi},  $\Gamma (f) \le \Gamma (f_3) +1 \le n$.
This completes the proof.
\end{proof}

\begin{rem}\label{rem:duality}
\begin{enumerate}
{\em \item
For a collection $\E$ of spaces,
the $\E\dash$Clapp-Puppe category of $f$ fits nicely into
our axiomatic framework.   Precisely,  let
$\cal{T}(\E)$  be the collection  of all maps  $g:X\to Y$
such that $g \simeq vu$, where $u:X\to E$ and $v:E\to Y$, for some $E
\in \E$.
It can be shown that $\E\dash\cat_{\mathrm{CP}}(f) = \Gamma
_{\cal{T}(\E)}(f)$ for all maps $f$.
\item
The axioms of Definition \ref{defn:axi} can be dualized
in the sense of Eckmann-Hilton.  This just consists of replacing Axiom
3 with an appropriate homotopy pull-back axiom.  As in Definition
\ref{defn:maxi}, we set $\Lambda _{\SSS}$ equal to the maximum
all functions which satisfy the dual axioms.  If $\E $ is any collection
of spaces, we   define
$$
\SSS ^*(\E) =
\left\{
g: X\to Y\, \left|
\begin{array}{l}
\exists\ \mathrm{maps}\ h_1: Y\to E\ \mathrm{and}\ h_2:E\to   X
\\
\mathrm{with}\ E\in \E\ \mathrm{and}\ h_2h_1 g \simeq \id_X
\end{array}
\right.\right\}.
$$

Then $\Lambda _{\SSS ^*(\E)}(f)$,
which we can denote by $\E\dash$cosecat$(f)$, is the dual of the
$\E\dash$sectional category of $f$.  It would be interesting to
investigate this invariant. }
\end{enumerate}
\end{rem}

\section{Inequalities}\label{section:inequalities}

In this section we consider the questions
raised in the introduction dealing with the $\E$-sectional category of a
composition and the relationship between the various invariants.  The
main theorem of this section is Theorem
\ref{thm:FH}.

\subsection{Composition of Maps}

The composition results in this subsection are
proved by elementary covering arguments.

\begin{prop}\label{prop:obvious}
For any maps $f:X\to Y$ and $g:Y\to Z$,
$$
\E\dash\sct(g) \le \E\dash\sct(gf).
$$
\end{prop}

The proof is obvious, and hence omitted.

\medskip

The next result deals with maps that have sections.

\begin{prop} If $f:X\to Y$ and $g:Y\to Z$,
then
\begin{enumerate}
\item   If  $f$ has a section, then $\E\dash\sct (gf)  =  \E\dash\sct
(g)$.
\item   If  $g$ has a section, then $\E\dash\sct (gf) \leq \E\dash\sct
(f)$.
\end{enumerate}
\end{prop}

\begin{proof}
First assume that $f$ has a section $s$.  Let $U\sseq Z$ be
$\E$-section-categorical with respect to $g$.  We claim that $U$ is
also $\E$-section-categorical with respect to $gf$.  This
follows
easily from the following homotopy-commutative diagram
$$
\xymatrix{
&& X \ar[d]^f
\\
E  \ar[rru]^{s v} \ar[rr]^v    && Y \ar[d]^g
\\
U\ar[rr]\ar[u]^u && Z.}
$$
This proves
$\E\dash\sct (gf)  \le  \E\dash\sct (g)$, and
Proposition \ref{prop:obvious} completes the proof of (1).
When $g$ has a section $t$, one proves (2)
by showing that if $U\sseq Y$ is
$\E$-section-categorical
with respect to $f$, then $t^{-1}(U)$
is $\E$-section-categorical with respect to $gf$.
\end{proof}

Next we  derive additional
results on the sectional category of a composition.
We begin with a lemma.

\begin{lem}\label{lem:mc}
If $A\maprt{} Y \maprt{j}Y\cup CA$ is a mapping cone sequence
and $f:X\maprt {} Y$ is a map, then
$$
\E\dash\sct (jf)  \le \E\dash\sct (f) +1.
$$
\end{lem}

\begin{proof} Consider the diagram
$$
\xymatrix{
       \ast  \ar[d]^{ } &&  \ast  \ar[ll]_{ } \ar[rr]^{ } \ar[d]^{ } && X
\ar[d]^{f}
\\
       \ast  && A \ar[ll]_{ }\ar[rr]^{ } && Y
\\  }
$$
which induces the map $jf:X\to Y\cup CA$ of homotopy pushouts.
The result now follows from Theorem \ref{thm:po}.
\end{proof}

Our main result gives an upper bound for the $\E$-sectional category of a
composition.

\begin{thm}\label{thm:FH}
If $f:X\maprt{} Y$ and $g:Y\maprt{}  Z$, then
$$
\E\dash\sct (gf)  \le {\rm cat}_{\mathrm{FH}}(g) + \E\dash\sct (f) .
$$
\end{thm}

\begin{proof}
Suppose $\cat_{\mathrm{FH}}(g) = n$.  Then we can choose a
decomposition of $g$ \cite[\S3]{ASS}
$$
\shortlabeldecomp{Y}{g}{Z}{n}{K}{i}{j}{s_n}
$$
where (1) $g \simeq g_nj_{n-1}\cdots j_1 j_0$,
(2) $s_ng \simeq j_{n-1}\cdots j_1 j_0$,
(3) $g_ns_n \simeq {\mathrm{id}}$ and
(4) there is a mapping cone sequence
$
\xymatrix@1{
A_i   \ar[r]        &
Y_i    \ar[r]^-{j_i} & Y_{i+1}
}
$.
Since $g$ is dominated by $j_{n-1}\cdots j_1 j_0$,
it follows that $ gf$ is dominated
by $j_{n-1}\cdots j_1 j_0f$.  By Lemma \ref{lem:mc},
\begin{eqnarray*}
\E{\rm -seccat}(gf)
& \le &\E\dash\sct ( j_{n-1}\cdots j_1 j_0f) \\
& \le & \E\dash\sct ( j_{n-2}\cdots j_1 j_0f)+1 \\
      & \vdots \\ &
\le & \E\dash\sct (f)  + n.
\end{eqnarray*}
\end{proof}

\medskip

We conclude this subsection by showing that
the inequality in Theorem \ref{thm:FH} can be strict.

\begin{ex}
{\em
Let $f:X\to Y$ be any map with
$\E\dash\sct (f) >0$, let $Z =*$ and let $g:Y\to Z$ be the
constant map.  Then
$$
\E\dash\sct (gf) = {\rm cat}(*) = 0 < \cat_{\mathrm{FH}}(g) +
\E\dash\sct (f).
$$
}
\end{ex}

\subsection{Relations Between Invariants}

We begin this subsection with an inequality which provides
a lower bound for the Clapp-Puppe category of a composition.

\begin{prop}\label{prop:CPcomp}
If $f:X\to Y$ and $g:Y\to Z$, then
$$
\E\dash\cat_{\mathrm{CP}} (g) + 1
\leq
(\sct (f)+ 1)(\E\dash\cat_{\mathrm{CP}}(gf) + 1).
$$
\end{prop}

\begin{proof}
Write $\sct(f) = n$ and $\E\dash\cat_{\mathrm{CP}}(gf)= m$.
Then we have a cover $\{ U_0, U_1,\ldots ,U_n\} $ of $Y$
and sections $s_i: U_i\to X$ for each $i$. Now
$$
\E\dash\cat_{\mathrm{CP}}( g|_{U_i} ) = \E\dash\cat_{\mathrm{CP}}
( gfs_i   )
\leq
\E\dash\cat_{\mathrm{CP}}(gf) = m,
$$
so $U_i = V_{i0}\cup \cdots \cup V_{im}$ with each $g|V_{ij}$
factoring through some member of $\E$.  Then $Y = \{ V_{ij}\, | \,
i = 0, 1, 2, \dots ,n; \, \; j = 0, 1,2,\ldots , m \}$
is desired cover of $Y$ by  $(n+1)(m+1)$ $\E$-categorical subsets.
\end{proof}

If $*$ is a constant map, then clearly $\E\dash\cat_{\mathrm{CP}}(*)=0$.
Thus if $\E = \{ *\} $ and $gf \simeq *$, then Proposition
\ref{prop:CPcomp}
yields the inequality
$$
\cat _{\mathrm{BG}}(g) \le \sct (f)
$$
of Berstein-Ganea \cite[Prop 2.6]{BG}.  More generally,
if $gf \simeq *$, then
$$
\E\dash\cat_{\mathrm{CP}}(g) \le \sct (f)
$$
for any collection $\E$ by Proposition \ref{prop:CPcomp}.  This gives
$\E\dash\cat_{\mathrm{CP}}(g) \le\E\dash\sct (f) $ if $gf \simeq *$.
However this inequality holds without the latter assumption.

\begin{prop}\label{prop:compo}
Let $f:X\to Y$ and $g:Y\to Z$.  Then
$$
\E\dash\cat_{\mathrm{CP}} (g)\leq \E\dash\cat_{\mathrm{CP}} (Y)  \leq
\E\dash\sct(f).
$$
\end{prop}

\begin{proof}
The left inequality follows immediately.
The right inequality follows from the observation that if $\{U_0,U_1, 
\ldots, U_n\}$
is an open cover of $Y$ which is
$\E$-section-categorical for $f$, then it is an open cover for
the $\E\dash$Clapp-Puppe category of id$_Y$.
\end{proof}

Next we turn to a corollary to Theorem \ref{thm:FH}.

\begin{cor}\label{cor:scc}
For any map $f:X\to Y$,
$$
\E\dash\sct(f) \leq  \cat_{\mathrm{FH}}(f) +
\E\dash\cat_{\mathrm{CP}}(X)  .
$$
In particular,  $\sct(f) \leq  \cat_{\mathrm{FH}}(f)$.
\end{cor}

\begin{proof}
Apply Theorem \ref{thm:FH} to the composition
$\xymatrix@1{X\ar[r]^{\id} & X \ar[r]^f & Y}$
and use Proposition \ref{prop:id}(1).
\end{proof}

This corollary is used to prove the following relationships
between the invariants $\E\dash\sct$, $\E\dash\cat_{\mathrm{CP}}$
and $\E\dash\cat_{\mathrm{FH}}$.

\begin{prop}\label{prop:in}
Let $\E $ be a collection and $f:X\to Y$ and $g:Y\to Z$ be maps.  Then
\begin{enumerate}
\item
$\E\dash\cat_{\mathrm{CP}}(g) \le \cat_{\mathrm{FH}}(f) +
\E\dash\cat_{\mathrm{CP}} (X)$.
\item
If $X\in \E$, then
$\E\dash\cat_{\mathrm{CP}} (g) \le \cat_{\mathrm{FH}}(f)$.
\item
If $gf \simeq *$, then
$\cat_{\mathrm{BG}}(g) \le \cat_{\mathrm{FH}}(f)$.
\end{enumerate}
\end{prop}

\begin{proof}
By Proposition \ref{prop:compo} and Corollary \ref{cor:scc}, we have
$$
\begin{array}{rcl}
\E\dash\cat_{\mathrm{CP}}(g)  & \leq   & \E\dash\sct(f)
\\
&\leq &
\cat_{\mathrm{FH}}(f)  +  \E\dash\cat_{\mathrm{CP}}(X).
\end{array}
$$
Also (2) follows from (1) since
$\E\dash\cat_{\mathrm{CP}}(X) = 0$.
For (3) we specialize to $\E =\{ \mathrm{all\; spaces}\} $ in Corollary
\ref{cor:scc}
and use the Berstein-Ganea inequality $\cat_{\mathrm{BG}}(g) \le \sct
(f)$
mentioned above.
\end{proof}

In the special case $\E= \{ *\}$, Proposition \ref{prop:in}(1) reduces
to
$$
\cat_{\mathrm{BG}} (g) \leq \cat_{\mathrm{FH}}(f) + \cat(X).
$$

Next we relate sectional category and the Clapp-Puppe category.

\begin{prop}\label{prop:CP}
For any collection $\E$ and any map $f:X\to Y$,
$$
\E\dash\cat_{\rm CP}(f) \le \E\dash\sct (f).
$$
\end{prop}

\begin{proof}
Let $\E\dash\sct (f) =n$ and let  $\{U_0,U_1, \ldots ,U_n\}$ be a
$\E$-section-categorical cover of
$Y$.  Then there are spaces
$E_i \in \E$ and maps $u_i:U_i\to E_i$ and $v_i:E_i \to X$
such that $fv_iu_i \simeq j_i$, where $j_i:U_i\inclds Y$ is the
inclusion map.  We set $V_i = f^{-1}(U_i)$.  Then $\{V_0,V_1, \ldots
,V_n\}$
is an open cover of $X$ and the maps
$$
\xymatrix{
V_i\ar[r]^{f|_{V_i}} &
U_i\ar[r]^{u_i} &
E_i\ar[r]^{v_i} &
X\ar[r]^{f} &
Y
}
$$
show that $\E\dash\cat_{\rm CP}(f) \le n.$
\end{proof}

\begin{rem}
{\em
This proposition does {\it not} show that $\cat_{\rm BG}(f)
\linebreak \le \sct(f)$.
This is because $\E\dash\cat_{\rm CP}(f) =
    \cat_{\mathrm{BG}}(f)$ when
$\E = \{ * \}$ and
$\E\dash\sct(f) = \sct(f)$ when $\E =  \{ \mathrm{all\ spaces}\}$.
In a sense, $\E\dash\cat_{\rm CP}$ and $\E\dash\sct$ are dual to each
other.
}
\end{rem}

\vskip 1in

\noindent  Dartmouth College

\noindent Hanover, NH 03755

\noindent Martin.Arkowitz@Dartmouth.edu

\vskip .3in

\noindent Western Michigan University

\noindent Kalamazoo, MI 49008

\noindent Jeff.Strom@wmich.edu

\end{document}